\documentclass[12pt]{article}

\usepackage{epsf}
\usepackage{amsmath}
\usepackage{amsthm}
\usepackage{amssymb}

\renewcommand{\epsfsize}[2]{\textwidth}

\theoremstyle{plain}
\newtheorem{thm}{Theorem}[section]
\newtheorem{prop}[thm]{Proposition}

\newtheorem{lem}[thm]{Lemma}

\newtheorem{intlemnp}[thm]{Lemma}
\newenvironment{lemnp}{\begin{intlemnp}}{\qed \end{intlemnp}}
\newtheorem{intthmnp}[thm]{Theorem}
\newenvironment{thmnp}{\begin{intthmnp}}{\qed \end{intthmnp}}

\theoremstyle{definition}

\theoremstyle{remark}
\newtheorem{rem}[thm]{Remark}

\title{Hyperbolic automorphisms of free groups}
\date{}
\author{Peter Brinkmann}

\begin{document}

\maketitle

\begin{abstract}
We prove that an automorphism $\phi:F\rightarrow F$ of a finitely generated
free group $F$ is hyperbolic in the sense of Gromov if it has no nontrivial
periodic conjugacy classes. This result was previously claimed (but not proved)
in \cite{combi}.
\end{abstract}

\section{Introduction}
Let $F$ be a finitely generated free group. We fix a basis once and for all,
and denote by $|.|$ the word length with respect to this basis.
An automorphism $\phi:F\rightarrow F$ is
said to be {\it hyperbolic in the sense of Gromov} (or just
{\it hyperbolic}) if there exist numbers $M>0$ and $\lambda>1$ such that
$$ \lambda |g| \leq \max \{ |\phi^M(g)|,|\phi^{-M}(g)| \} $$
for all $g\in F$.

An automorphism $\phi:F\rightarrow F$ is called {\it atoroidal} if it has
no nontrivial periodic conjugacy classes. This definition is motivated by
the fact that the mapping torus $F\rtimes_\phi \mathbb Z$ of such an
automorphism contains no subgroups isomorphic to
$\mathbb{Z}\oplus\mathbb{Z}$.

The following theorem is the main result of this paper.
It was previously claimed in \cite{combi}, and a proof in
the special case of irreducible automorphisms appeared in \cite{bfgafa},
\footnote{In \cite{bfgafa2}, the authors state that the main results of
\cite{bfgafa} are only proved for
irreducible automorphisms with irreducible powers, but a
close inspection of their proofs
shows that this additional hypothesis is not necessary for the results
quoted in this paper.} but the general question remained open until now.

\begin{thm}\label{main}
If $\phi:F\rightarrow F$ is an atoroidal automorphism,
then $\phi$ is hyperbolic.
\end{thm}

Theorem \ref{main} completes the proof of the following theorem.
\begin{thm}
Let $\phi$ be an automorphism of $F$. Then the following statements
are equivalent:
\begin{enumerate}
\item \label{eq1} The mapping torus $F\rtimes_\phi \mathbb Z$ of $\phi$ is
hyperbolic.
\item \label{eq2} $\phi$ is hyperbolic.
\item \label{eq3} $\phi$ is atoroidal.
\end{enumerate}
\end{thm}

The equivalence between (\ref{eq1}) and (\ref{eq2}) was established in
\cite{combi,combi2}
and \cite{isoper}. (\ref{eq2}) $\Rightarrow$ (\ref{eq3}) follows
immediately from the definitions, and (\ref{eq3}) $\Rightarrow$ (\ref{eq2})
is precisely our Theorem \ref{main}.

The train track techniques developed in \cite{hb1,tits1} will be our most
important tool, and in Section \ref{train}, we review the results needed in 
this paper.  The key feature of our proof of Theorem \ref{main}
is that we study the growth of the length of paths under iterates of
$\phi^{\pm 1}$ by analyzing only one train track map
representing some iterate of $\phi$, rather than two train track maps
representing $\phi$ and $\phi^{-1}$. This was originally suggested 
by Martin Lustig, and it simplifies the approach considerably
as the relationship between train track maps representing $\phi$ and
train track maps representing $\phi^{-1}$ remains somewhat mysterious.

I would like to thank Mladen Bestvina, Steve Gersten and Martin Lustig for
many helpful discussions. I am indebted to the referee for several
helpful suggestions, and I would also like to express my gratitude to Sean
Sather-Wagstaff.

\section{Hyperbolic automorphisms and homotopy equivalences of graphs}

In this section, we list some basic lemmas that will allow us to use
train track techniques in the proof of Theorem \ref{main}. We also introduce
some notation that will be useful later.

\begin{lemnp}[\cite{bfgafa}]\label{chyp}
Let $\phi:F\rightarrow F$ be an atoroidal automorphism.
Then $\phi$ is hyperbolic if there exist numbers
$M>0$ and $\lambda>1$ such that
$$ \lambda ||g|| \leq \max \{ ||\phi^M(g)||,||\phi^{-M}(g)|| \} $$
for all $g\in F$, where $||.||$ denotes the length
of the conjugacy class of $g$, i.\ e., the length of the shortest word in
the conjugacy class of $g$.
\end{lemnp}

\begin{rem}
The lemma shows that hyperbolicity is in fact a property of the
outer automorphism represented by $\phi$, which will allow us to prove
Theorem \ref{main} by geometric methods. Lemma \ref{bilip} rephrases
the problem in terms of geometry.
\end{rem}

Let $G$ be a finite graph with fundamental group $F$.
We will always assume that a homotopy equivalence $f:G\rightarrow G$
maps vertices to vertices, and we will only consider
graphs whose vertices have valence at least two.
Following the conventions of \cite{tits1}, we will refer to a map
$\rho:[0,1]\rightarrow G$ as a {\it path} if it is either constant or an
immersion, and we reserve the word {\it circuit}
for immersions $\sigma:S^1\rightarrow G$.

We do not require paths to start or end at vertices.
For a path or circuit $\rho$ in $G$, let $[f(\rho)]$ denote
the path or circuit homotopic (relative end points if $\rho$ is a path) to
the composition of $\rho$ with $f$.
For a subpath $\rho$ of a path or circuit $\sigma$, let $[f^k(\rho)]_\sigma$
denote the maximal subpath of $[f^k(\rho)]$ contained in $[f^k(\sigma)]$.

Given a homotopy equivalence $f:G\rightarrow G$ of a graph and a path or
circuit $\rho$ in $G$, we
denote by $\rho^{-k}$ a path or circuit in $G$ with the property that
$[f^k(\rho^{-k})]=\rho$. Such a path $\rho^{-k}$ always exists, but it may not
be unique if
$\rho$ is not a circuit. However, all the statements we will prove will be
independent of the choice of $\rho^{-k}$.
Given some metric on $G$,
we denote the length of $\rho$ by $L(\rho)$, and we restrict our attention
to homotopy equivalences that map edges to paths of positive length.

\begin{lem}\label{bilip}
Let $f:G\rightarrow G$ be a homotopy equivalence representing an
outer automorphism $\mathcal O\in Out\ F$.
If there exist numbers $M$ and $\lambda>1$ such that
$$\lambda L(\sigma)\leq \max\{L([f^M(\sigma)]),L(\sigma^{-M})\}$$
holds for all circuits $\sigma$ in $G$, then $\mathcal O$ is
hyperbolic.
\end{lem}

\begin{proof} There exists a constant $C>0$ such that
$$C^{-1}L(\sigma)\leq ||x|| \leq CL(\sigma),$$
for all conjugacy classes $x$ in $F$ and circuits $\sigma$ in $G$ representing
$x$. Choose $K$ such that $\lambda^K>C$. We conclude that
$${\frac{\lambda^K}{C}}||x|| \leq \max\{||{\mathcal O}^{KM}(x)||,
||{\mathcal O}^{-KM}(x)||\}$$
for all conjugacy classes $x$.
\end{proof}

In light of Lemma \ref{bilip}, we call a homotopy equivalence
{\it hyperbolic} if it represents a hyperbolic automorphism.

\section{Train tracks}\label{train}

In this section, we review the theory of train tracks developed in
\cite{hb1,tits1}. We will restrict our attention to the collection of
those results that will be used in this paper.

Oftentimes, a homotopy equivalence $f:G\rightarrow G$ will respect a
{\it filtration}
of $G$, i.\ e., there exist subgraphs $G_0=\emptyset\subset G_1 \subset \cdots
\subset G_k=G$ such that for each filtration element $G_r$,
the restriction of $f$ to $G_r$
is a homotopy equivalence of $G_r$. The subgraph
$H_r=\overline{G_r\setminus G_{r-1}}$ 
is called the {\it $r$-th stratum} of the filtration.
We say that a path $\rho$ has {\it nontrivial intersection} with a
stratum $H_r$ if $\rho$ crosses at least one edge in $H_r$.

If $E_1,\cdots,E_m$ is the collection of edges in some stratum $H_r$, the
{\it transition matrix} of $H_r$ is the nonnegative
$m\times m$-matrix $M_r$ whose $ij$-th
entry is the number of times the $f$-image of $E_j$ crosses $E_i$, regardless
of orientation. $M_r$ is said to be {\it irreducible} if
for every tuple $1\leq i,j \leq m$, there exists some exponent $n>0$ such that
the $ij$-th entry of $M_r^n$ is nonzero.
If $M_r$ is irreducible, then it has a maximal real eigenvalue
$\lambda_r\geq 1$
(see \cite{seneta}).  We call $\lambda_r$ the
{\it growth rate} of $H_r$.

Given a homotopy equivalence $f:G\rightarrow G$, we can always find a
filtration of $G$ such that each transition matrix is either a zero matrix
or irreducible. A stratum $H_r$ in such a filtration is called
{\it zero stratum} if $M_r=0$. $H_r$ is called {\it exponentially growing}
if $M_r$ is irreducible with $\lambda_r>1$, and it is called
{\it polynomially growing} if $M_r$ is irreducible with $\lambda_r=1$.

An unordered pair of edges in $G$ originating from the same vertex is called
a {\it turn}. A turn is called {\it degenerate} if the two edges are equal.
We define a map
$Df:\{\text{turns in } G\}\rightarrow \{\text{turns in } G\}$ by sending each
edge in a turn to the first edge in its image under $f$. A turn is called
{\it illegal} if its image under some iterate of $Df$ is degenerate, {\it legal}
otherwise.

An edge path $\alpha=E_1E_2\cdots E_s$ is said to contain the turns
$(\bar{E_i},E_{i+1})$ for $1\leq i <s$. $\alpha$ is said to be legal if
all its turns are legal, and a path $\alpha\subset G_r$ is $r$-legal if
no illegal turn in $\alpha$ involves an edge in $H_r$.

A path $\rho$ in $G$ is said to be a {\it (periodic) Nielsen path} if
$\rho$ is not constant and if
$[f^k(\rho)]=\rho$ for some $k>0$. The {\it period} of $\rho$ is the smallest
such exponent $k$.
A {\it pre-Nielsen} path is a path
whose image under some iterate of $f$ is a Nielsen path.
A Nielsen path is called
{\it indivisible} if it cannot be written as a concatenation of two Nielsen
paths.

A decomposition of a path or circuit $\sigma=\sigma_1\cdots\sigma_s$ into
subpaths
is called a {\it splitting} if $[f^k(\sigma_i)]=[f^k(\sigma_i)]_\sigma$ for
all $k, i$.

The following theorem was proved in \cite{hb1}.
\begin{thmnp}[{\cite[Theorem 5.12]{hb1}}]
Every outer automorphism $\mathcal O$ of $F$ is represented by a
homotopy equivalence $f:G\rightarrow G$ such that each exponentially
growing stratum $H_r$ has the following properties:
\begin{enumerate}
\item If $E$ is an edge in $H_r$, then the
first and last edges in $f(E)$ are contained in $H_r$.
\item If $\beta$ is a nontrivial path in $G_{r-1}$ with endpoints in
$G_{r-1}\cap H_r$, then $[f(\beta)]$ is nontrivial.
\item If $\rho$ is an $r$-legal path, then $[f(\rho)]$ is an $r$-legal path.
\end{enumerate}
\end{thmnp}
We call $f$ a {\it relative train track map}.

An outer automorphism $\mathcal O$ of $F$ is called {\it reducible} if it
preserves the conjugacy class of a proper free factor of $F$. $\mathcal O$
is called {\it irreducible} if it is not reducible. If $\mathcal O$ is
irreducible, then it has a relative train track representative
$f:G\rightarrow G$ whose filtration has only one nonempty element $H_1=G$,
with irreducible transition matrix. The properties of relative train tracks
show that for every edge $E$ of $G$, the image $f^n(E)$ is an immersion
for all $n>0$. In this case, we call $f$ a {\it train track map} (or 
absolute train track map), and we denote the growth rate of $H_1=G$ by
$\lambda$.

We now construct a metric on $G$. If $H_r$ is an exponentially growing
stratum, then its transition matrix $M_r$ has a unique
positive left eigenvector $v_r$
(corresponding to $\lambda_r$) whose smallest entry equals one
(see \cite{seneta}).
For an edge $E_i$ in $H_r$, the eigenvector $v_r$ has an entry $l_i>0$
corresponding to $E_i$. We choose a metric on $G$ such that $E_i$ is
isometric to an interval of length $l_i$, and such that edges in zero strata
or in polynomially growing strata are isometric to an interval of length one.
Note that if $\rho$ is a path whose endpoints are vertices, then the number
of edges in $\rho$ provides a lower bound for $L(\rho)$. Moreover, if $f$
is an absolute train track map, then $f$ expands the length of legal paths
by the factor $\lambda$.

For our purposes, the properties of relative train track maps are not strong
enough, so we will use the notion of improved train track maps
constructed in \cite{tits1}. We only list the properties used in this
paper.

\begin{thmnp}[{\cite[Theorem 5.1.5]{tits1}}]\label{ttimproved}
For every outer automorphism $\mathcal O$ of $F$, there exists
an exponent $k>0$ such that $\mathcal O^k$ is represented by a relative
train track map $f:G\rightarrow G$ with the following additional properties:
\begin{enumerate}
\item Every periodic Nielsen path has period one.
\item $H_r$ is a zero stratum if and only if it is the union of the
contractible components of $G_r$.
\item If $H_{r-1}$ is a zero stratum, then $H_r$ is an exponentially growing
stratum.
\item \label{polyim} If $H_r$ is a polynomially growing stratum,
then $H_r$ consists of
a single edge $E_r$, and $f(E_r)=E_ru_r$, where $u_r\subset G_{r-1}$.
\item \label{oneINP} If $H_r$ is an exponentially growing stratum, then
there is at most
one indivisible Nielsen path in $G_r$ that intersects $H_r$ nontrivially.
\end{enumerate}
\end{thmnp}
We call $f$ an {\it improved} relative train track map.

The following lemma is an immediate consequence of \cite[Lemma 5.1.7]{tits1}.
\begin{lemnp}\label{npprop}
Suppose that $f:G\rightarrow G$ is an improved train track map representing
an atoroidal automorphism, with an exponentially growing stratum $H_r$.
If $\rho\subset G_r$ is a Nielsen path whose first and last edges are
contained in $H_r$, then the endpoints of $\rho$ are distinct, and if
both endpoints are contained in $G_{r-1}$, then at least one of them is
contained in a contractible component of $G_{r-1}$.
\end{lemnp}

The following lemma will turn out to be crucial in the proof of Theorem
\ref{main}.
\begin{lem}\label{legalrest}
Let $f:G\rightarrow G$ be an improved train track map representing an
atoroidal outer automorphism.
If $H_r$ is an exponentially growing stratum
and if $\rho$ is an indivisible Nielsen path in $G_r$ with nontrivial
intersection with $H_r$, then the endpoints of $\rho$ are distinct and
at least one of them is not contained in $H_r\cap G_{r-1}$.
\end{lem}

\begin{proof}
Since $\rho$ is indivisible, its initial edge and its terminal edge are
contained in $H_r$.
Suppose that both endpoints of $\rho$ are contained in $H_r\cap G_{r-1}$.
By Lemma \ref{npprop}, one of them is contained in a contractible
component of $G_{r-1}$. Let $v$ denote this endpoint. By Theorem
\ref{ttimproved}, $H_{r-1}$ is necessarily a zero stratum, and we have
$v\in H_{r-1}$.
This implies that $f(v)\in G_{r-2}$. Since $H_{r-1}$ is the collection
of contractible components of $G_{r-1}$, we conclude that
$H_{r-1}\cap G_{r-2}=\emptyset$, which implies that $v\neq f(v)$.
This contradicts our assumption that $\rho$ is a Nielsen path.
\end{proof}

Finally, we state a lemma from \cite{tits1} that simplifies the study of
paths intersecting strata of polynomial growth. If $H_r=\{E_r\}$ is a
polynomially growing stratum, then {\it basic paths}
of height $r$ are of the form $E_r\gamma$,
$\gamma \bar{E_r}$, or $E_r\gamma \bar{E_r}$, where $\gamma$ is a path in
$G_{r-1}$ with endpoints in $H_r$.

\begin{lemnp}[{\cite[lemma 4.1.4]{tits1}}]\label{pgsplit}
Let $f:G\rightarrow G$ be an improved train track map with a polynomially
growing stratum $H_r$. If $\sigma$ is a circuit in $G_r$, then
it splits as a concatenation of basic paths of height $r$ and paths in
$G_{r-1}$.
\end{lemnp}

\begin{rem}\label{howtosplit}
In fact, part \ref{polyim} of Theorem \ref{ttimproved} implies that subdividing
$\sigma$ at the initial endpoints of all occurrences of $E_r$ and at the
terminal endpoints of all occurrences of $\bar{E}_r$ yields a
splitting of $\sigma$ into basic paths of height $r$ and paths in $G_{r-1}$.
\end{rem}

\section{Bounded cancellation}

Thurston's bounded cancellation lemma is one of the fundamental tools in this
paper. We state it in terms of homotopy equivalences of graphs.

\begin{lemnp}[Bounded cancellation lemma, see \cite{cooper}]\label{bcc}
Let $f:G\rightarrow G$ be a homotopy equivalence. There exists a constant
$C_f$, depending only on $f$, with the property that for any path $\rho$ in $G$
obtained by concatenating two paths $\alpha, \beta$, we have
$$ L([f(\rho)])\geq L([f(\alpha)]) + L([f(\beta)]) - C_f. $$
\end{lemnp}

Bounded cancellation allows us to draw conclusions
about the growth of sufficiently long paths under iterates of hyperbolic
homotopy equivalences. We make this precise in the following lemma.

\begin{lem}\label{pathhyp}
Let $f:G\rightarrow G$ be a homotopy equivalence, and let $G'$ be a subgraph
of $G$ such that the restriction of $f$ to $G'$ is a hyperbolic homotopy
equivalence of $G'$.

Then there exist constants $L^C, \lambda>1$ and $N$ with the property that,
if $\rho\subset G'$ is a subpath of some circuit $\sigma$ in $G$ and if
the length of $\rho$ is at least $L^C$, we have
$$ \lambda L(\rho)\leq \max \{L([f^N(\rho)]_\sigma), L(\alpha)\},$$
where $\alpha$ is some subpath of $\sigma^{-N}$ satisfying
$[f^N(\alpha)]_{(\sigma^{-N})}=\rho$.
\end{lem}

\begin{proof}

Let $g:G\rightarrow G$ be a homotopy inverse of $f$.
Since the restriction of $f$ to $G'$ is hyperbolic,
there exist numbers $\lambda', N$ such that
$$\lambda' L(\sigma')\leq \max \{L([f^N(\sigma')]),L([g^N(\sigma')])\}$$
for all circuits $\sigma'$ in $G'$.
\par

There exists a path $\tau\subset G'$ such that $\sigma'=[\rho\tau]$ is an
immersed circuit in $G'$, satisfying $L(\tau)\leq 2 \text{diam}(G')$
and $L(\sigma')\geq L(\rho)-\text{diam}(G')$.
We can find some constant $K_1$,
depending only on $G$ and $f$, such that $L([g^N(\tau)])\leq K_1$,
$L([f^N(\tau)])\leq K_1$ and $\lambda'\text{diam}(G')\leq K_1$.
We distinguish two cases that are not mutually exclusive.
\begin{enumerate}
\item $\lambda' L(\sigma')\leq L([f^N(\sigma')])$. In this case, the bounded
cancellation lemma tells us that
\begin{equation*}
\begin{split}
\lambda'L(\rho) \leq \lambda'(L(\sigma') +\text{diam}(G'))
\leq L([f^N(\sigma')]) + K_1 \\
\leq L([f^N(\rho)]_{\sigma'})+2K_1 
\leq L([f^N(\rho)]_\sigma)+2K_1+2C_{f^N},
\end{split}
\end{equation*}
where $C_{f^N}$ is the bounded cancellation constant of $f^N$.

\item $\lambda' L(\sigma')\leq L([g^N(\sigma')])$. The same reasoning as in
the previous case shows that
$\lambda' L(\rho)\leq L([g^N(\rho)]_\sigma)+2K_1+2C_{g^N}$, where
$C_{g^N}$ is the bounded cancellation constant of $g^N$.
Since $f^N$ and $g^N$ are homotopy inverses of each other, we can find some
constant $K_2$ such that $|L(\beta^{-N})-L([g^N(\beta)])|\leq K_2$
holds for all
paths $\beta$ and preimages $\beta^{-N}$ ($\beta^{-N}$ is a preimage of
$\beta$ under $f^N$). We conclude that
$$ \lambda'L(\rho)\leq L(\alpha)+2K_1+K_2+2C_{g^N}, $$
where $\alpha$ is some subpath of $\sigma^{-N}$ satisfying
$[f^N(\alpha)]_{(\sigma^{-N})}=\rho$.
\end{enumerate}

Let $K=2K_1+\max\{2C_{f^N}, K_2+2C_{g^N}\}$.
Choose $L^C$ large enough to satisfy $\lambda' L^C - K > L^C$ and
let $\lambda=\frac{\lambda' L^C - K}{L^C}$. Clearly, $\lambda>1$.
If the length of $\rho$ is at least $L^C$, we conclude that
$$\lambda L(\rho) \leq \lambda' L(\rho) - K
\leq \max \{L([f^N(\rho)]_\sigma), L(\alpha)\}. $$
\end{proof}

We call $L^C$ the {\it critical length} of the triple $(f,G,G')$.

For train track maps, there is a related
concept of critical length. Let $f:G\rightarrow G$ be a train
track map with growth rate $\lambda$ and bounded cancellation constant
$C_f$. If $\beta$ is a legal path in $G$ whose length satisfies
$\lambda L(\beta)-2C_f>L(\beta)$ and $\alpha, \gamma$ are paths such
that the concatenation $\alpha\beta\gamma$ is locally injective, then
the length of the segment in $[f^n(\alpha\beta\gamma)]$ corresponding to
$\beta$ will tend to infinity as $n$ tends to infinity. In this situation, the
critical length is the infimum of the lengths satisfying the above inequality,
i.\ e., $\frac{2C_f}{\lambda-1}$.

\section{The irreducible case}
Throughout this section, let $f:G\rightarrow G$ denote a
train track map representing an exponentially growing irreducible
outer automorphism
$\mathcal O$ of $F$, with growth rate $\lambda>1$.
We equip $G$ with the metric constructed in Section \ref{train}.

For a path or circuit
$\rho$ in $G$, let $\mathcal{L}(\rho)$ denote the length
of the longest legal segment of $\rho$ (recall that $L(\rho)$ denotes the
length of $\rho$).
Let $i(\rho)$ denote the number of illegal turns in $\rho$.

As in the previous section, we denote by $\rho^{-k}$ a path or circuit in $G$
with the property that $[f^k(\rho^{-k})]=\rho$.

Let $C_f$ denote the bounded cancellation constant of $f$, and let
$L^c={\frac{2C_f} {\lambda-1}}$ be the critical length of $f$.

We will use the following lemma from \cite{bfgafa}.
\begin{lemnp}\label{tricho1}
For all $L>0$ there exists
an exponent $M>0$ such that if $\rho$ is any path in $G$, one of the following
holds:
\begin{enumerate}
\item $[f^M(\rho)]$ has a legal segment of length greater than $L$.
\item $[f^M(\rho)]$ has fewer illegal turns than $\rho$.
\item $\rho$ can be expressed as a concatenation $\tau_1\rho'\tau_2$, where
$L(\tau_1)\leq 2L$, $L(\tau_2)\leq 2L$, $i(\tau_1)\leq 1$,
$i(\tau_2)\leq 1$, and
$\rho'$ splits as a concatenation of pre-Nielsen paths with one illegal
turn each.
\end{enumerate} 
\end{lemnp}

In order to study the length of preimages $\rho^{-k}$ of a path $\rho$,
we will need an upper bound on $\mathcal{L}({\rho^{-k}})$ in terms of
$\mathcal{L}(\rho)$.
We make this precise in the following lemma.

\begin{lem}\label{bw1}
For all paths or circuits $\rho$ in $G$ and exponents $k>0$, we have
$$\mathcal{L}({\rho^{-k}})<{\frac{\mathcal{L}(\rho)} {\lambda^k}}+L^c.$$ 
In particular, this implies that $\mathcal{L}({\rho^{-k}})<
\mathcal{L}(\rho)+L^c$.
\end{lem}

\begin{proof} We will show by induction that
$$ \mathcal{L}({\rho^{-k}})\leq {\frac{\mathcal{L}(\rho)
 +2C_f(1+\lambda+\cdots+\lambda^{k-1})}
{\lambda^k}}.$$
For $k=1$, the bounded cancellation lemma implies that
$$\lambda \mathcal{L}({\rho^{-1}}) - 2C_f\leq \mathcal{L}(\rho)
\quad \Leftrightarrow \quad
\mathcal{L}({\rho^{-1}})\leq {\frac{\mathcal{L}(\rho)+2C_f}{\lambda}},$$
so the claim holds for $k=1$.

Assume that the claim is true for some $k$. Again, the bounded cancellation
lemma tells us that
$$
\mathcal{L}({\rho^{-(k+1)}})\leq
\frac{\mathcal{L}({\rho^{-k}})+2C_f}{\lambda} \leq
{\frac{\mathcal{L}(\rho) +2C_f(1+\lambda+\cdots+\lambda^{k})}{\lambda^{k+1}}}
$$
by induction. Hence, we conclude that

$$\mathcal{L}({\rho^{-k}})\leq
{\frac{\mathcal{L}(\rho)+2C_f{\frac{\lambda^k-1}{\lambda-1}}}
{\lambda^k}}={\frac{\mathcal{L}(\rho)}{\lambda^k}}+2C_f{\frac{1-\lambda^{-k}}
{\lambda-1}}<{\frac{\mathcal{L}(\rho)}{\lambda^k}}+{\frac{2C_f}{\lambda-1}}.$$
\end{proof}

While we measure the growth of paths under forward iteration by means of the
path metric in $G$, our measure of growth under backward iteration will be
the number of illegal turns. We make this precise in Lemmas \ref{backgrowth}
and \ref{illen}. Lemma \ref{backgrowth} is a stronger version of
\cite[Lemma 2]{lustig}.

\begin{lem}\label{backgrowth}
Let $f:G\rightarrow G$ be an improved train track map inducing an atoroidal
outer automorphism $\mathcal O$.
Given some number $L_0>0$, there exists some exponent $M>0$, depending only on
$L_0$, such that for any path $\rho$ with $\mathcal{L}(\rho)\leq L_0$ and
$i(\rho)\geq 4$,
we have
$$ \left(\frac{8}{7}\right)^n i(\rho)\leq i(\rho^{-nM}) $$
for all $n>0$.
\end{lem}

\begin{proof}
Given $L_0$, choose an exponent $M$ according to Lemma \ref{tricho1}, for
$L=L_0+L^c$. 
Express $\rho$ as a concatenation of paths $\rho_1,\cdots,\rho_s,\tau$ such
that $i(\rho_i)=4$ and $i(\tau)<4$. There exist preimages $\rho^{-M}_1,
\cdots, \rho^{-M}_s,\tau^{-M}$ such that $\rho^{-M}$ is their concatenation.

We claim that $i(\rho^{-M}_i)\geq 5$ for all $i$. Suppose otherwise,
i.\ e., $i(\rho^{-M}_i)=4$ for some $i$. Because
of our choice of $M$, Lemmas \ref{tricho1} and \ref{bw1} imply
that $\rho^{-M}_i$ can be written as a
concatenation $\tau_1\rho'\tau_2$,
where $\rho'$ splits as a concatenation of two
pre-Nielsen paths with one illegal turn each.
This implies that
for some exponent $k\geq 0$, $[f^k(\rho')]$
contains the concatenation of
two Nielsen paths, which is impossible
because of Lemma \ref{legalrest} and Theorem \ref{ttimproved}, part
\ref{oneINP}.

Hence, $i(\rho^{-M})\geq 5s + i(\tau) \geq
\left(\frac{8}{7}\right) i(\rho)$, and the lemma follows by induction.
\end{proof}

The following lemma establishes an elementary connection between the length of
a path and the number of illegal turns.
\begin{lemnp}\label{illen}
Given some $L>0$, there exists some constant $C>0$ such that for all
paths $\rho$ with $1\leq \mathcal{L}(\rho)\leq L$ and $i(\rho)>0$, we have
$$C^{-1}i(\rho)\leq L(\rho) \leq C\ i(\rho).$$
\end{lemnp}

A version of the
following special case of Theorem \ref{main} has already been proved in
\cite{bfgafa}. We present a new proof.
\begin{thm}\label{main1}
If $f:G\rightarrow G$ is an improved train track map representing an
irreducible, atoroidal outer automorphism of $F$, then $f$ is hyperbolic.
\end{thm}

\begin{proof}
Fix some $L_0>L^c$.

Let $\sigma$ be a nontrivial circuit in $G$.
We will distinguish several cases, and in each case we will show that there
exist numbers $N>0$, $\lambda>1$ and $\epsilon>0$ such that there exists
a collection $S$ of subpaths of $\sigma$ with the following properties:
\begin{enumerate}
\item For every integer $n>0$ and for every $\rho\in S$, we have
$$\lambda^n L(\rho)\leq \max\{[f^{nN}(\rho)]_\sigma,L(\alpha)\},$$
where $\alpha$ is a subpath of $\sigma^{-nN}$ such that
$[f^{nN}(\alpha)]_{(\sigma^{-nN})}=\rho$. We say that $\rho$ has the
{\it desired growth}.
\item There is no overlap between distinct paths in $S$.
\item The sum of the lengths of the paths in $S$ is at least
$\epsilon L(\sigma)$.
\end{enumerate}
If the numbers $N$, $\lambda$ and $\epsilon$ depend only on the case in
question, but not on $\sigma$, then the theorem follows immediately because the
growth of the subpaths in $S$ provides a lower bound for the growth of
$\sigma$.

We distinguish the following cases.
\begin{enumerate}
\item $i(\sigma)=0$ or $\frac{L(\sigma)}{i(\sigma)}\geq L_0$.
If $i(\sigma)=0$,
then $\sigma$ is legal, so it has the desired growth under forward
iteration. Otherwise, let $S$ be the collection of maximal legal
subpaths of $\sigma$
of length at least $L_0$. The choice of $L_0$ and Lemma \ref{bcc}
guarantee that the subpaths in $S$ have the desired growth under
forward iteration, so we only have to show that they account for a definite
fraction of the length of $\sigma$. An elementary computation will verify this.

Let $l$ be the length of the longest path whose endpoints are vertices
and whose length is strictly less than $L_0$. If $L(S)$ denotes the sum of
the lengths of the segments in $S$, we have $L(\sigma)-L(S)\leq i(\sigma)l$
and $L(\sigma)\geq i(\sigma)L_0$. This implies that
$$ \frac{L(S)}{L(\sigma)}=1-\frac{L(\sigma)-L(S)}{L(\sigma)}\geq
1-\frac{l}{L_0},$$
independently of $\sigma$.

\item $\frac{L(\sigma)}{i(\sigma)} < L_0$. There are two subcases to
consider.
\begin{enumerate}
\item $i(\sigma)\geq 4$. In this case, we define $S'$ to be the set of
subpaths left after removing from $\sigma$ the maximal legal segments of
length greater than $6L_0$. Then we obtain $S$ by removing from $S'$ the
subpaths with fewer than four illegal turns.

Lemma \ref{backgrowth} and Lemma \ref{illen}
(with $L=6L_0$) show that the elements of $S$ have the desired growth under
backward iteration, so we only have to show that $S$ accounts for a definite
positive fraction of the length of $\sigma$. This fraction is minimal if
$S$ contains only one subpath ($S$ cannot be empty) and if all the paths in
$S'\setminus S$ have exactly three illegal turns. We first find a lower bound
for the number of legal segments of the path in $S$: Let $a$ be the number
of legal segments of $\sigma$ of length greater than $6L_0$, and let $b$
be the number of legal segments in $S$. Then $a\leq \frac{i(\sigma)}{6}$
and $i(\sigma)=b+a+4(a-1)$, which implies $b\geq \frac{i(\sigma)}{6}$.
Moreover, the number of edges in a path provides a lower bound for the length
of that path, so $b$ is also a lower bound for $L(S)$,
and we conclude that
$$ \frac{L(S)}{L(\sigma)}\geq \frac{b}{L(\sigma)} \geq
\frac{i(\sigma)}{6L(\sigma)} \geq \frac{1}{6L_0}.$$

\item $i(\sigma) < 4$. In this case, the length of $\sigma$ is bounded by
$3L_0$, so there are only finitely many circuits to consider. Since
$\mathcal O$ is atoroidal, the length of all circuits tends
to infinity under forward iteration, and we can easily find an exponent $N$
with the property that, say,  $L([f^n(\sigma)])\geq 4L_0$ for all circuits
of length at most $3L_0$ and for all $n\geq N$.

\end{enumerate}
\end{enumerate}

The cases considered above account for every circuit $\sigma$.
This completes the proof.
\end{proof}

\section{The reducible case --- exponentially growing strata}
Throughout this section, let $f:G\rightarrow G$ be a relative train track map
representing an outer automorphism $\mathcal O\in Out\ F$.
The notation used in this section will be consistent with the notation in the
previous section; the subscript $r$ will indicate the stratum of $G$ under
consideration.

If $H_r$ is an exponentially growing stratum, let $\lambda_r$ be the
corresponding growth rate. We equip $G$ with the metric constructed in
Section \ref{train}.

Let $\rho$ be a path or circuit in $G_r$.
Following \cite{tits1}, we
denote by $\rho\cap H_r$ the ordered sequence of oriented edges of $H_r$
crossed by $\rho$.
We will refer to the total length of $\rho\cap H_r$ as the
$r$-{\it length} of $\rho$, denoted by
$L_r(\rho)$. Similarly, $i_r(\rho)$ denotes the number of $r$-illegal
turns in $\rho$, and $\mathcal{L}_r(\rho)$ stands for the $r$-length of
the \mbox{($r$-)}longest \mbox{$r$-}legal segment of $\rho$.
Let $L^c_r=\frac{2C_f}{\lambda_r-1}$ be the critical $r$-length,
where $C_f$ is the bounded cancellation constant of $f$.
The relative train track property implies that $f$ expands the $r$-length
of $r$-legal paths by the factor $\lambda_r$.

As in the previous section, we denote by $\rho^{-k}$ a path or circuit in $G$
with the property that $[f^k(\rho^{-k})]=\rho$.

The following lemma is a straight-forward generalization of Lemma \ref{bw1}.
\begin{lemnp}\label{bw2}
Let $f:G\rightarrow G$ be a train track map with an exponentially
growing stratum $H_r$, and let $\rho$ be a path or circuit in $G_r$. Then
$$ \mathcal{L}_r({\rho^{-k}}) < \mathcal{L}_r(\rho) + L^c_r. $$
\end{lemnp}

The following generalization of Lemma \ref{tricho1} is the main technical
result of this section, and it will be crucial for our analysis of backward
growth in the reducible case.
\begin{prop}\label{tricho2}
Let $f:G\rightarrow G$ be a relative train track map, and
let $H_r$ be an exponentially growing stratum. For all $L>0$, there
exists some exponent $M>0$ such that if $\rho$ is a path in $G_r$ with
$L_r(\rho)\geq 1$,
one of the following three statements holds:
\begin{enumerate}
\item $[f^M(\rho)]$ has an $r$-legal segment of $r$-length greater than $L$.
\item $[f^M(\rho)]$ has fewer $r$-illegal turns than $\rho$.
\item $\rho$ can be expressed as a concatenation $\tau_1\rho'\tau_2$, where
$L_r(\tau_1)\leq 2L$, $L_r(\tau_2)\leq 2L$, $i_r(\tau_1)\leq 1$,
$i_r(\tau_2)\leq 1$, and
$\rho'$ splits as a concatenation of pre-Nielsen paths (with one
$r$-illegal turn each) and segments in $G_{r-1}$. \label{case3}
\end{enumerate}
\end{prop}

In order to prove Proposition \ref{tricho2}, we will need the
following version of a well-known fact from Ramsey theory.
\begin{lem}\label{ramsey}
For all natural numbers $K, N_0, Q$ there exists some $M$ such that for all
maps $f:\{1,\cdots,M\}\rightarrow \{1,\cdots, K\}$ there exist numbers
$n$ and $N\geq N_0$ such that $f(n)=f(n+N)=\cdots=f(n+QN)$.
\end{lem}

\begin{proof} By \cite[page 55, Theorem 2]{ramsey},
there exists some number $M'$ such that
for all $f:\{1,\cdots,M'\}\rightarrow \{1,\cdots, K\}$ there exist numbers
$n$ and $N\geq 1$ such that $f(n)=f(n+N)=\cdots=f(n+QN)$. Now $M=M'N_0$ has
the desired property.
\end{proof}

\begin{proof}[Proof of Proposition \ref{tricho2}]
Fix some $N_0$ such that $\lambda_r^{N_0}\geq L$.
Although there may be infinitely many paths $\rho'\subset G_r$
whose endpoints are
vertices and whose $r$-length is at most $3L$, there is only a finite number
$K$ of intersections $\rho'\cap H_r$ of such paths with $H_r$. 

Choose $M$ according to Lemma \ref{ramsey}, with $Q=4$ and $K, N_0$ as above.
We will show that $M$ is the desired exponent.

Let $\rho$ be a path in $G_r$ with $L_r(\rho)\geq 1$.
Suppose that the first
two statements do not hold for $[f^M(\rho)]$. We want to show that 
the third statement is satisfied.
In order to avoid case  distinctions, we
assume that $i_r(\rho)\geq 4$; the proof in the case $i_r(\rho)<4$ is a
straight-forward modification of the following argument.

Let $\rho^k_1,\cdots,\rho^k_m$ denote the $r$-legal segments of $[f^k(\rho)]$
for $k=0,\cdots,M$. By assumption, the $r$-length of each $\rho^k_i$ is
bounded by $L$. Moreover, as the turn between two consecutive subpaths is
$r$-illegal, the last edge in $\rho^k_i$ is contained in
$H_r$ if $i<m$, and the first edge in $\rho^k_i$ is contained in $H_r$
if $i>1$.

Fix some $1<i<m-2$.\footnote{We do not consider the segments $\rho_1$ and
$\rho_m$ in the following argument because the initial endpoint of $\rho_1$
(resp.\ the terminal endpoint of $\rho_m$) may not be a vertex. The paths
$\tau_1$ and $\tau_2$ in the third statement account for $\rho_1$ and
$\rho_m$.}
By Lemma \ref{ramsey}, there
exist numbers $n$ and $N\geq N_0$ such that
$\rho^n_i\cap H_r=\rho^{n+N}_i\cap H_r =\cdots=\rho^{n+4N}_i\cap H_r$,
$\rho^n_{i+1}\cap H_r=\rho^{n+N}_{i+1}\cap H_r =\cdots
=\rho^{n+4N}_{i+1}\cap H_r$, and
$\rho^n_{i+2}\cap H_r=\rho^{n+N}_{i+2}\cap H_r =\cdots
=\rho^{n+4N}_{i+2}\cap H_r$.

An elementary topological argument shows that there exist subpaths $\alpha'$
of $\rho^n_i\rho^n_{i+1}\rho^n_{i+2}$ with the property
$\alpha'\cap H_r$=$[f^N(\alpha')]\cap H_r$ and $i_r(\alpha')=2$.
Let $\alpha$ be the shortest such subpath.
We will show that $\alpha$ can be expressed as a
concatenation $\alpha_1\gamma\alpha_2$, where $\alpha_1, \alpha_2$ are
pre-Nielsen paths, and $\gamma$ is a path that is constant or
contained in $G_{r-1}$.

There exists a unique shortest subpath $\alpha_1$ of $\alpha$ such that
$\alpha_1$ contains the first illegal turn of $\alpha$ and 
$\alpha_1\cap H_r=[f^N(\alpha_1)]\cap H_r$. Similarly, let $\alpha_2$ be the
shortest subpath of $\alpha$ such that $\alpha_2$ contains the second illegal
turn of $\alpha$ and $\alpha_2\cap H_r=[f^N(\alpha_2)]\cap H_r$.
Note that the extremal (i.e., initial and terminal)
edges of $\alpha_1, \alpha_2$ are (possibly partial) edges in $H_r$.
We have $\alpha=\alpha_1\gamma\alpha_2$ for some path $\gamma\subset G_r$.
If $\gamma$ were a path of positive $r$-length, this would imply
that $L_r(f^N(\gamma))>L_r(\gamma)$, contradicting our choice of $\alpha,
\alpha_1$ and $\alpha_2$. We conclude that $\gamma$ is contained in $G_{r-1}$
or constant.

We claim that
$\alpha_1\cap H_r=[f^N(\alpha_1)]\cap H_r=\cdots=
[f^{4N}(\alpha_1)]\cap H_r$.
In order to see this, we need to understand
the cancellation that occurs between
the two maximal $r$-legal subpaths of $f^N(\alpha_1)$.
In the tightening process, the
terminal edge of the first subpath cancels with the initial edge of the second
subpath until the last edge of the first subpath
forms a nondegenerate turn with the first edge of the second subpath.
Since $[f^N(\alpha_1)]$ contains an $r$-illegal turn,
the resulting turn is necessarily $r$-illegal; in particular, it is contained
in $H_r$.

This shows that the part of $f^N(\alpha_1)$ that is cancelled 
is completely determined by $\alpha_1\cap H_r$. Similarly, the part of
$f^N([f^{N}(\alpha_1)])$ that is cancelled is completely determined by
$[f^{N}(\alpha_1)]\cap H_r$, etc. Since
$\alpha_1\cap H_r=[f^N(\alpha_1)]\cap H_r$,
this shows that $\alpha_1\cap H_r=[f^{2N}(\alpha_1)]\cap H_r$. We
conclude that $\alpha_1\cap H_r=[f^N(\alpha_1)]\cap H_r=\cdots=
[f^{4N}(\alpha_1)]\cap H_r$.

Let $E$ denote the first (possibly partial) edge of $\alpha_1$.
The map $f^N$ expands the $r$-length of $E$ by $\lambda_r^N$,
and it maps
vertices to vertices, so $f^N(E)$ contains at least one entire edge in $H_r$,
which implies that $f^{2N}(E)$ has $r$-length at least $L$.
The same argument applies to
the last edge of $\alpha_1$, which implies that $[f^{2N}(\alpha_1)]$ is
completely determined by the extremal edges of $\alpha_1$. Applying this
argument to $[f^{2N}(\alpha_1)]$ and $[f^{4N}(\alpha_1)]$, we conclude that
$[f^{2N}(\alpha_1)]=[f^{4N}(\alpha_1)]$, hence
$[f^{2N}(\alpha)]$ is a Nielsen path.
The same argument shows that $[f^{2N}(\alpha_2)]$ is a Nielsen path.

Repeating this argument for all indices $1<i<m-2$, we conclude that $\rho$
splits as a concatenation
$\rho=\tau_1\beta_1\gamma_1\beta_2\gamma_2\cdots\beta_{m-3}\tau_2$, where
$\tau_1$ and $\tau_2$ are as in the third statement, the paths $\beta_i$ are
pre-Nielsen, and the paths $\gamma_i$ are contained in $G_{r-1}$ or constant.
This completes the proof.
\end{proof}

We will need the following relative versions of Lemmas \ref{backgrowth}
and \ref{illen}.

\begin{lem}\label{bgrowth2}
Assume that $\mathcal O$ is atoroidal and that
$f:G\rightarrow G$ is an improved train track representative with an
exponentially growing stratum $H_r$.
Given some number $L_0>0$, there exists some exponent $M>0$, depending only on
$L_0$ and $H_r$,
such that for any path $\rho$ in $G_r$
with $\mathcal{L}_r(\rho)\leq L_0$ and $i_r(\rho)\geq 5$, we have
$$ \left(\frac{10}{9}\right)^n i_r(\rho)\leq i_r(\rho^{-nM}) $$
for all $n>0$.
\end{lem}

\begin{proof}
Given $L$, choose an exponent $M$ according to Lemma \ref{tricho2},
with $L=L_0+L^c_r$.
Express $\rho$ as a concatenation of paths $\rho_1,\cdots,\rho_s,\tau$ such
that $i_r(\rho_i)=5$ and $i_r(\tau)<5$. There exist preimages of $\rho^{-M}_1,
\cdots, \rho^{-M}_s,\tau^{-M}$ such that $\rho^{-M}$ is their concatenation.

We claim that $i_r(\rho^{-M}_i)\geq 6$ for all $i$. Suppose otherwise, i.\ e.,
$i_r(\rho^{-M}_i)=5$ for some $i$.
Because of
our choice of $M$, Lemma \ref{tricho2} implies that $\rho^{-M}_i$ can be
written as a concatenation $\tau_1\rho'\tau_2$, and
$\rho'$ splits as a concatenation of three pre-Nielsen paths (with one
illegal turn each) with (possibly empty) segments in $G_{r-1}$ in between. 
This implies that for some exponent $k\geq 0$, $[f^k(\rho')]$ contains
three Nielsen paths with segments in $G_{r-1}$ in between, which is impossible
because of Lemma \ref{legalrest} and Theorem \ref{ttimproved}, part
\ref{oneINP}.

Hence, $i_r(\rho^{-M})\geq 6s+i_r(\tau)\geq \frac{10}{9} i_r(\rho)$, and the
lemma follows by induction.
\end{proof}

\begin{lemnp}\label{illen2}
Suppose $H_r$ is an exponentially growing stratum.
Given some $L>0$, there exists some constant $C>0$ such that for all
paths $\rho\subset G_r$ with $1\leq \mathcal{L}_r(\rho)\leq L$ and
$i_r(\rho)>0$, we have
$$C^{-1}i_r(\rho)\leq L_r(\rho) \leq C\ i_r(\rho).$$
\end{lemnp}

\section{Proof of the main theorem}

The following proposition and Lemma \ref{bilip} immediately imply Theorem
\ref{main}.

\begin{prop}
If $f:G\rightarrow G$ is an improved relative train track map representing an
atoroidal outer automorphism, then $f$ is hyperbolic.
\end{prop}

\begin{proof}
We will proceed by induction up through the filtration of $G$ (as in the
previous sections, we equip $G$ with the metric constructed in Section
\ref{train}). The
restriction of $f$ to $G_1=H_1$ is a homotopy equivalence, and we claim
that $H_1$ is of exponential growth. If it were a zero stratum, this would
imply that $f(H_1)\subset G_0$, but $G_0=\emptyset$.
If it were of polynomial growth, it would give rise to a
nontrivial fixed conjugacy class. We conclude that $H_1$
is of exponential growth, so this initial case follows from
Theorem \ref{main1}.

Now assume that the restriction of $f$ to $G_{r-1}$ is hyperbolic.
We choose constants $L^C, \lambda$ and $N$ according to Lemma
\ref{pathhyp} for the triple $(f,G_r,G_{r-1})$.

As in the proof of Theorem \ref{main1}, we will distinguish several cases,
and in each case we will find a collection $S$ of subpaths having the desired
growth and accounting for a definite positive fraction of the length of the
circuit in question.

For the inductive step, we distinguish three main cases, depending on the
stratum $H_r$.

\begin{enumerate}

\item $H_r$ is a zero stratum. Then $H_r$ is the collection of
contractible components of $G_r$ (see Theorem \ref{ttimproved}).
This implies that any nontrivial circuit in $G_r$ is
contained in $G_{r-1}$, so there is nothing to show in this case.

\item $H_r$ is an exponentially growing stratum.
We fix some length $L_0>\max \{L^C, L^c_r\}$.
Let $\sigma$ be a circuit in $G_r$ with nontrivial intersection with $H_r$.
If $H_{r-1}$ is a zero stratum, we let $\sigma_1=\sigma\cap G_{r-2}$,
$\sigma_0=\sigma\cap H_{r-1}$ and $\sigma_2=\sigma\cap H_r$. Otherwise,
let $\sigma_1=\sigma\cap G_{r-1}, \sigma_0=\emptyset$
and $\sigma_2=\sigma\cap H_r$.

If $H_{r-1}$ is a zero stratum, then it is the collection of contractible
components of $G_{r-1}$ (see Theorem \ref{ttimproved}). Consider a subpath
$\rho$ of $\sigma$ that is contained in $H_{r-1}$. If $\rho$ is maximal,
i.\ e., if $\rho$ is not a proper subpath of another subpath of $\sigma$
that is contained in $H_{r-1}$, then the edges preceding and following
$\rho$ in the edge circuit $\sigma$ are contained in $H_r$.

This implies that
$\frac{L(\sigma_0)}{L(\sigma_2)}\leq \text{diam}(G)$
no matter whether $H_{r-1}$ is a zero stratum or not. Hence,
$L(\sigma)\leq L(\sigma_1)+(1+\text{diam}(G))L(\sigma_2)$.

As in the proof of Theorem \ref{main1},
we will decompose $\sigma$ into subpaths whose
growth we understand. We consider several cases.

\begin{enumerate}
\item \label{case1}
$\frac{L(\sigma_1)}{L(\sigma_2)}\geq L_0$. In this case, there will be segments
of length at least $L_0$ in $\sigma_1$, and the inductive hypothesis and
Lemma \ref{pathhyp} show that they have the desired growth.
Hence, it suffices to show that those segments 
account for some definite fraction $\epsilon > 0 $ of the length of
$\sigma$, where $\epsilon$ does not depend on the choice of $\sigma$.
An elementary computation will verify this.

Let $A$ be the total
length of all segments of length at least $L_0$ in $\sigma_1$, $B$ the total
length of the remaining segments in $\sigma_1$, and let $C=L(\sigma_2)$.
Then our assumption implies $A+B\geq L_0C$. Moreover, if $m$ denotes the number
of segments in $\sigma_1$, we have $m\leq C$ and $\frac{A+B}{m}\geq
\frac{L_0C}{m}\geq L_0$.

We want to find a lower bound for $\frac{A}{L(\sigma)}$. Using the inequalities
derived so far, we conclude that
$$\frac{A}{L(\sigma)}\geq
\frac{A}{A+B+(1+\text{diam}(G))C}\geq
\frac{AL_0}{(A+B)(1+L_0+\text{diam}(G))}.$$
Hence, we only need to find a lower bound for $\frac{A}{A+B}$.
Let $l$ be the length of the longest path whose endpoints are vertices
and whose length is strictly less than $L_0$.
Then $B\leq ml$ and $A+B\geq mL_0$, and we conclude that
$\frac{A}{A+B}=1-\frac{B}{A+B}\geq 1-\frac{l}{L_0},$
independently of $\sigma$.

\item $\frac{L(\sigma_1)}{L(\sigma_2)} < L_0$. In this case, significant growth
will occur in $H_r$, and as in the proof of
Theorem \ref{main1}, we distinguish two subcases depending on whether forward
or backward growth dominates.
\begin{enumerate}
\item $\frac{L_r(\sigma)}{i_r(\sigma)}\geq L_0$. In analogy with case
\ref{case1}, we only need to show that $r$-legal segments of $r$-length at
least $L_0$ account for a definite fraction of the length of $\sigma$,
which can be accomplished with a computation very similar to the one in
case \ref{case1}.

\item $\frac{L_r(\sigma)}{i_r(\sigma)} < L_0$. As in the proof of Theorem
\ref{main1}, we consider two subcases.
\begin{enumerate}
\item $i_r(\sigma)\geq 5$. We define $S'$ to be the set of subpaths left after
removing from $\sigma$ the maximal $r$-legal subpaths of $r$-length greater
than $7L_0$, and we obtain $S$ from $S'$ by removing subpaths with fewer
than five $r$-illegal turns.

Lemma \ref{bw2} (with $L=7L_0$) and Lemmas \ref{bgrowth2} and \ref{illen2}
show that the paths in $S$ have the desired growth under backward
iteration. An argument very similar to the one in the proof of Theorem
\ref{main1} shows that the sum of the lengths of the paths in $S$ accounts
for a definite positive fraction of the length of $\sigma$,
so we are done in this case.
\item $i_r(\sigma)<5$. Only finitely many circuits $\sigma$ fall into this
category, and the same argument as in the proof of Theorem \ref{main1}
shows that
they have the desired growth under forward iteration.

\end{enumerate}
\end{enumerate}

\end{enumerate}

\item $H_r$ is a polynomially growing stratum. Recall (see Section \ref{train})
that $H_r$ contains
only one edge $E_r$, and that basic paths of height $r$ are of the form
$E_r\gamma$, $\gamma \bar{E_r}$, or $E_r\gamma \bar{E_r}$, where $\gamma$
is a path in $G_{r-1}$ with endpoints in $H_r$.

We fix some $L_0>L^C$.
Let $\sigma$ be a circuit in
$G_r$ with nontrivial intersection with $H_r$. Using Lemma \ref{pgsplit} and
Remark \ref{howtosplit}, we obtain a splitting of $\sigma$ by subdividing
$\sigma$ at the initial endpoints of all occurrences of $E_r$ and at the
terminal endpoints of all occurrences of $\bar{E}_r$. The subpaths of $\sigma$
obtained in this way are either basic paths of height $r$ or paths in
$G_{r-1}$, and the endpoints of all subpaths are contained in $H_r$.

We first show that all basic paths of height $r$ have the desired growth under
sufficiently high iterates of $f^{\pm 1}$. Let $\rho$ be a basic path of height
$r$. Since a basic path of the form $\gamma \bar{E_r}$ can be turned into
a basic path of the form $E_r\gamma$ by reversing its orientation,
we only have to distinguish two cases.

\begin{enumerate}
\item $\rho=E_r\gamma\bar{E_r}$ with $\gamma\subset G_{r-1}$. If
$L(\gamma)\geq L_0$, the inductive hypothesis and Lemma \ref{pathhyp} prove
the claim, so it suffices to consider the case $L(\gamma)<L_0$. The endpoints
of $\gamma$ are equal, and we denote by $\tau$ the circuit defined by $\gamma$.
In general, $\tau$ may be shorter than $\gamma$ because initial and terminal
edges of $\gamma$ may cancel. However, we have $L(\tau)\geq 1$ and
$L(\gamma)-L(\tau)<L_0$. Moreover, the growth of $\tau$ under iterates of $f$
provides a lower bound for the growth of $\rho$ under iterates of $f$, so the
inductive hypothesis proves the claim in this case.

\item $\rho=E_r\gamma$ with $\gamma\subset G_{r-1}$. As in the previous case,
we may assume that $L(\gamma)<L_0$. We first show that $[f(\rho)]\neq \rho$.
Suppose otherwise. Then the endpoints of $\rho$ cannot be equal because
$\mathcal O$ is atoroidal. This implies that the endpoints of
$E_r$ are distinct and $\gamma$ starts and ends at the terminal endpoint of
$E_r$. However, this is impossible as it implies that $\rho$ is a basic
path of the form $E_r\gamma\bar{E_r}$ (see Lemma \ref{pgsplit} and Remark
\ref{howtosplit}).

We conclude that $\lim_{n\rightarrow \infty}L([f^n(\rho)])=\infty$. As
there are only finitely many paths of length less than $L_0$, we conclude that
the circuits in this category have the desired growth under forward
iteration.
\end{enumerate}

We have shown that basic paths of height $r$ have the desired growth,
as do paths in $G_{r-1}$ if their length is
at least $L_0$. This leaves us with those subpaths in the splitting of $\sigma$
that are contained in $G_{r-1}$ and whose length is less than $L_0$, but we
can safely disregard them because there are at least as many basic paths of
height $r$ as there are subpaths in $G_{r-1}$. This completes the proof.

\end{enumerate}

\end{proof}

\bibliographystyle{is-alpha}
\bibliography{all}
\bigskip
{\sc Department of Mathematics, University of Utah\\
Salt Lake City, UT 84112, USA\\}
\\
{\it E-mail:} brinkman\@@math.utah.edu

\end{document}